\documentclass[12pt]{amsart}
\usepackage{amsmath}
\usepackage{amsthm}
\usepackage{amsfonts}
\usepackage{graphicx}
\usepackage{amssymb}
\newcommand{\arcangle}{\mathord{<\mspace{-9mu}\mathrel{)}\mspace{2mu}}}

\begin{document}

\title[Minimal Surfaces from Rigid Motions]
{Minimal Surfaces from Rigid Motions}
\author[Jens Hoppe ]{JENS HOPPE}
\address{Braunschweig University, Germany}
\email{jens.r.hoppe@gmail.com}

\begin{abstract}
Equations are derived for the shape of a hypersurface in $\mathbb{R}^N$
for which a rigid motion yields a minimal surface in $\mathbb{R}^{N+1}$.
Some elementary, but unconventional, aspects of the classical case
$N=2$ (solved by H.F. Scherk in 1835) are discussed in some detail.
\end{abstract}

\maketitle
\noindent
The Ansatz
\begin{equation}\label{eq1}
x(t,\varphi) = \begin{pmatrix}
R(f(t))\vec{u}(\varphi) \\ t
\end{pmatrix},
\end{equation}
with $R(f) = e^{f\left( \begin{smallmatrix}
0 &-1 \\ 1 & 0
\end{smallmatrix}\right)}$
a t(ime)-dependent rotation gives
\begin{equation}\label{eq2}
\begin{split}
\dot{x} & = \begin{pmatrix}
\dot{f}AR\vec{u} \\ 1
\end{pmatrix}, \;
x' = \begin{pmatrix}
R\vec{u}' \\ 0
\end{pmatrix} \\
(g_{ab}) & = \begin{pmatrix}
1 + \dot{f}^2 r^2 & \dot{f}(\vec{u} \times \vec{u}') \\ \cdot & \vec{u}'^2
\end{pmatrix},
\end{split}
\end{equation}
$g = det(g_{ab}) = \vec{u}'^2(1+ \dot{f}^2r^2 cos^2\phi)$, $r:= |\vec{u}|$, $\phi := \arcangle(\vec{u},\vec{u}')$,\\
and $\triangle := \dfrac{1}{\sqrt{g}}\partial_a(\sqrt{g}g^{ab}\partial_b)$ acting on $t$ the condition
\begin{equation}\label{eq3}
\ddot{f}|\vec{u}'|r^2 cos^2 \phi + (rS)' + \dot{f}^2 r^3 S' = 0
\end{equation}
where $S := sin \phi$; hence $\ddot{f} = 0$, as the special cases $(rS)' = 0$, resp. $r^3S' = 0$, resp. $\phi \equiv \frac{\pi}{2} = const$ can be excluded. So the planar curve described by $\vec{u}(\varphi)$ has to be rotated with constant angular velocity (we will take $f(t) = \omega t$ from now on, without loss of generality; the factor $\omega$ actually being somewhat superfluous, as multiplying (\ref{eq1}) by $\omega$ and redefining $t$ and $\vec{u}$ shows).
Integrating
\begin{equation}\label{eq4}
r'S + r(1+\omega^2 r^2)S' = 0
\end{equation}
gives
\begin{equation}\label{eq5}
\dfrac{\omega r sin \phi}{\sqrt{1+\omega^2 r^2}} = \varepsilon = const;
\end{equation}
instead of taking this constant of integration $(-1 < \varepsilon < +1)$ as the parameter, we will use the constant $\gamma_0 = \frac{\varepsilon}{\sqrt{1-\varepsilon^2}} \in \mathbb R$ that naturally arises from the original conservation law for $f  =\omega t$, resp. $\dot{g}_{ab} = 0$, i.e. rewriting $\triangle t = 0$ as
\begin{equation}\label{eq6}
0 = \dfrac{\partial}{\partial t}(\sqrt{g}g^{tt}) = -\dfrac{\partial}{\partial \varphi}(\sqrt{g}g^{21}) = \dfrac{\partial}{\partial \varphi}\underbrace{\left(\dfrac{\omega r sin \phi}{\sqrt{1+\omega^2 r^2 cos^2 \phi}}\right)}_{=: \gamma_0}  = 0,
\end{equation}
which (excluding $\gamma_0 = 0$, i.e. the helicoid) when written as
\begin{equation}\label{eq7}
\omega^2 r^2(\gamma^2 sin^2 \phi -1) = 1,
\end{equation}
$\gamma^2 = \frac{\gamma_0^2 +1}{\gamma_0^2} = \frac{1}{\varepsilon^2} > 1$ (cp.(\ref{eq5})) we will call the `shape-equation' (cp.\cite{1,2,3} for the relativistic case). Solving (\ref{eq7}) in arclength parametrization ($|\vec{u}\,'| = 1 \Rightarrow r cos \phi = \vec{u} \vec{u}' = \frac{1}{2}(r^2)' = rr'$, hence)
\begin{equation}\label{eq8}
cos \phi = r',
\end{equation}
starting $(s=0)$ at minimal distance from the origin ($\stackrel{\wedge}{=}$ maximal $\phi$), and going in
two different directions (i.e., contrary to the usual convention that the arclength parameter $s$ is necessarily $\geq 0$, having $s \in \mathbb R$), one gets
\begin{equation}\label{eq9}
\omega r = \pm \sqrt{\gamma_0^2 + \beta (\omega s)^2} \underset{\gamma_0 \neq 0}{=} \gamma_0 \sqrt{1 + \dfrac{(\omega s)^2}{\gamma_0^2(1+\gamma_0^2)}}
\end{equation}
$\beta := (\gamma_0^2 + 1)^{-1} \in (0,1]$, ($\gamma_0 = 0 , r(s) = |s|$, corresponding to the helicoid; noting that otherwise $\gamma_0$ and $\omega$ always have the same sign according to (\ref{eq6}), if $sin \phi \geqslant 0$; we will take the positive $s$ part of the infinite curve to be the one going against the clock). Instead of integrating (using the above convention)
\begin{equation}\label{eq10}
\theta' = \sqrt{\dfrac{1-r'^2}{r^2}} = |\omega|\dfrac{\sqrt{\gamma_0^2 + \beta(1-\beta)\omega^2 s^2}}{\gamma_0^2+\beta\omega^2 s^2} = \gamma_0 \dfrac{\omega \sqrt{1+\beta^2\omega^2 s^2}}{\gamma_0^2 + \beta \omega^2 s^2}
\end{equation}
to obtain
\begin{equation}\label{eq11}
\vec{u}(\varphi(s)) = r(s)\begin{pmatrix}
cos\theta(s) \\ sin\theta(s)
\end{pmatrix} = r \begin{pmatrix}
c \\ s
\end{pmatrix}
\end{equation}
one can also derive the curvature of the planar curve via
\begin{equation}\label{eq12}
\begin{split}
\vec{u}' & = r'\begin{pmatrix}
c \\ s
\end{pmatrix} + \theta' \begin{pmatrix}
-s \\ c
\end{pmatrix}\\
\vec{u}'' & = (r'' - r\theta'^2)\begin{pmatrix}
c \\ s
\end{pmatrix} + (r\theta'' + 2r'\theta')\begin{pmatrix}
-s \\ c
\end{pmatrix},
\end{split}
\end{equation}
true in \textit{any} parametrization.
In arclength parametrization ($\vec{u}' = \vec{c}\,' = \vec{e}_1,\, \vec{e_1}' = \kappa\vec{e}_2,\, \vec{e}_2 := -r\theta' \begin{pmatrix}
c \\ s
\end{pmatrix} +r'\begin{pmatrix}
-s \\ c
\end{pmatrix},\, det(\vec{e}_1, \vec{e}_2) = 1$)
\begin{equation}\label{eq13}
\kappa \underset{s>0}{=} \dfrac{r\theta'^2 - r''}{r\theta'} = \left(\dfrac{\sqrt{1-r'^2}}{r} - \dfrac{r''}{\sqrt{1-r'^2}}\right) = \dfrac{\kappa_0}{\sqrt{1+\mu^2s^2}},
\end{equation}
having used that for $r' \neq 0$ (\ref{eq7}), resp.
\begin{equation}\label{eq14}
\gamma^2(1-r'^2) = \dfrac{1}{\omega^2r^2} + 1
\end{equation}
implies
\begin{equation}\label{eq15}
\begin{split}
r'' = \dfrac{1}{\gamma^2 \omega^2 r^3}, \ 
r'^2 + r r'' = \beta,
rr' = s\beta;
\end{split}
\end{equation}
hence
\begin{equation}\label{eq16}
\kappa_0 = \gamma_0\omega\beta > 0,\; \mu = \beta \omega.
\end{equation}
With
\begin{equation}\label{eq17}
\sigma(s) := \int_0^s\kappa(v)dv = \dfrac{\kappa_0}{\mu} arcsinh (\mu s)
\end{equation}
the general solution of the Frenet-equations then gives
\begin{equation}\label{eq18}
\begin{split}
&\vec{u}(\varphi(s)) = \vec{c}(s) = \vec{c}(0) + \int_0^s {\begin{pmatrix}
 cos(\frac{\kappa_0}{\mu}arcsinh (\mu t)+ \sigma_0)\\
 sin(\frac{\kappa_0}{\mu}arcsinh (\mu t)+ \sigma_0)
 \end{pmatrix}}dt \\
 & = \vec{c}(0) + \int_0^{v = \frac{1}{\mu}arcsinh(\mu s)}{\begin{pmatrix}
 cosh(\mu u)cos(\kappa_0 u + \sigma_0)\\
 cosh(\mu u)sin(\kappa_0 u + \sigma_0)
 \end{pmatrix}}du\\
 & = \begin{pmatrix}
 \frac{\mu}{\mu^2+\kappa_0^2}sinh(\mu u)cos(\kappa_0 u + \sigma_0) + \frac{\kappa_0}{\mu^2+\kappa_0^2}cosh(\mu u)sin(\kappa_0 u + \sigma_0)\\
 \frac{\mu}{\mu^2+\kappa_0^2}sinh(\mu u)sin(\kappa_0 u + \sigma_0) - \frac{\kappa_0}{\mu^2+\kappa_0^2}cosh(\mu u)cos(\kappa_0 u + \sigma_0)
 \end{pmatrix}_0^v + \vec{c}(0)\\
 & = \dfrac{1}{\mu^2+\kappa_0^2} \begin{pmatrix}
 \mu sinh\, cos + \kappa_0 cosh\, sin \\
 \mu sinh\, sin - \kappa_0 cosh\, cos
 \end{pmatrix}(v) - \dfrac{\kappa_0}{\mu^2+\kappa_0^2} \begin{pmatrix}
 sin \sigma_0 \\ -cos \sigma_0
 \end{pmatrix} + \vec{c}(0)\\
 & = \vec{c}(0)+ \dfrac{1}{\mu^2+\kappa_0^2}\begin{pmatrix}
 cos & sin \\ sin & -cos
 \end{pmatrix} \begin{pmatrix}
 \mu sinh \\ \kappa_0 cosh
 \end{pmatrix} - \dfrac{\kappa_0}{\mu^2+\kappa_0^2}\begin{pmatrix}
 sin \sigma_0 \\ -cos \sigma_0
 \end{pmatrix}
\end{split}
\end{equation}
where $\mu t = sinh(\mu u)$, $\mu s = sinh(\mu v)$.
So
\begin{equation}\label{eq19}
\begin{split}
x(t, \varphi(s(v))) & = \begin{pmatrix}
R(\omega t) R(\sigma_0)[R(\kappa_0 v)\vec{w}(v) + \vec{w}_0] \\ t
\end{pmatrix}, \\
\vec{w}(v) & =  \dfrac{1}{\mu^2+\kappa_0^2} \begin{pmatrix}
\mu sinh(\mu v) \\ -\kappa_0 cosh(\mu v)
\end{pmatrix} = \begin{pmatrix}
a sinh(\tilde{v}) \\ b cosh (\tilde{v})
\end{pmatrix}\\
a\omega & = 1,\, b\omega = -\frac{\kappa_0}{\mu},\, \omega = \frac{\mu^2+\kappa_0^2}{\mu}
\end{split}
\end{equation}
where $\tilde{u} := \omega t + \kappa_0 v + \sigma_0, \, \tilde{v} = \mu v$ (i.e. $t = \frac{1}{\omega}(\tilde{u}-\frac{\kappa_0}{\mu}\tilde{v} - \sigma_0)$) so that for $\vec{w}_0 = \vec{0}$
\begin{center}
$\underset{^{\sim}}{x}(\underset{^{\sim}}{u},\underset{^{\sim}}{v}) = \begin{pmatrix}
R(\tilde{u})\vec{\underset{^{\sim}}{w}}(\tilde{v}) \\
a\tilde{u} + b\tilde{v} +c
\end{pmatrix},
$
\end{center}
i.e., dropping all $^{\sim}$'s, the constant $c$, and interchanging $u$ and $v$
\begin{equation}\label{eq20}
\begin{split}
x(u,v) & = \begin{pmatrix}
a sinh(u)cos(v) - b cosh(u)sin(v) \\
a sinh(u)sin(v) + b cosh(u)cos(v) \\
av + bu
\end{pmatrix}\\
a & = \dfrac{\mu}{\mu^2 + \kappa_0^2}, \,
b  = \dfrac{-\kappa_0}{\mu^2 + \kappa_0^2}.
\end{split}
\end{equation}
The curve $\vec{u}(\varphi)$ that is being rotated is
\begin{equation}\label{eq21}
\vec{r}(\tilde{v}) = \dfrac{1}{\omega}R(\gamma_0 \tilde{v})\begin{pmatrix}
sinh(\tilde{v}) \\ -\gamma_0cosh(\tilde{v})
\end{pmatrix};
\end{equation}
$\tilde{v} := arcsinh(\mu s)$, the picture on the next page showing $-\vec{r}$ for $\gamma_0 = \omega = -1$ (so $\tilde{v}$ and $s$ have opposite signs, as $\mu < 0$) resp., due to
\begin{equation}\label{eq22}
-\vec{r}_{-\gamma_0,-\omega}(\tilde{v}) = \vec{r}_{+\gamma_0,+\omega}(-\tilde{v}),
\end{equation}
also $-\vec{r}_{1,1}(-\tilde{v})$. As found by Bonnet \cite{4}, $u$ and $v$ are isothermal parameters for the minimal surfaces originally discovered by Scherk \cite{5}, $a=0$ corresponding to the Catenoid, $b=0$ to the helicoid;
\begin{equation}\label{eq23}
\begin{split}
x & = \begin{pmatrix}
R(v)\begin{pmatrix}
asinh(u)\\bcosh(u)
\end{pmatrix}\\
av+bu
\end{pmatrix}, \
x_u = \begin{pmatrix}
R\begin{pmatrix}
acosh\\bsinh
\end{pmatrix}\\
b
\end{pmatrix} = \begin{pmatrix}
R\vec{w}' \\b
\end{pmatrix}\\
x_v & =  \begin{pmatrix}
AR\vec{w} \\ a
\end{pmatrix}, \
x_u \cdot x_v = ab + \underbrace{\vec{w} \times \vec{w} '}_{\underset{ab(sinh^2 - cosh^2)}{\shortparallel}} = 0\\
x_v^2 & = \vec{w}^2 + a^2 = (a^2 + b^2)cosh^2 = \vec{w} '^2 +b^2 = x_u^2
\end{split}
\end{equation}
so that
\begin{equation}\label{eq24}
\partial_a \sqrt{g}g^{ab}\partial_b = \partial_u^2 + \partial_v^2
\end{equation}
(and the components of $x$ are obviously harmonic).
Going back to the original picture of a screw-motion applied to a planar curve i.e. $x = \big( \begin{smallmatrix} R(\omega z)\vec{r}(\eta) \\ z
\end{smallmatrix} \big)$ one finds for the Bonnet parametrization that $\vec{r}(\eta)$ (possible up to a fixed translation and rotation) is given by $\vec{r}(\eta) = R(\frac{-b}{a}\eta)  \big( \begin{smallmatrix} a sinh(\eta) \\ b cosh(\eta)
\end{smallmatrix} \big)$ resp.
\begin{equation}\label{eq25}
\omega\vec{r} = R(\gamma_0\eta)\begin{pmatrix}
sinh(\eta) \\ -\gamma_0 cosh(\eta)
\end{pmatrix} = \vec{v}(\eta) = \begin{pmatrix}
v_1 \\ v_2
\end{pmatrix}
\end{equation}
where $\eta \in \mathbb{R}$ $a,b \in \mathbb{R}$ $a\neq 0$, a one-parameter class of differently shaped curves
where $\gamma_0 \in \mathbb{R}$ is the constant arising in the original conservation law , $\gamma_0 = 0$ giving a straight line that is being rotated (the helicoid), while for $\gamma_0 \neq 0$ defining a curve starting at $(0, -\gamma_0)$, in each of the (positive/negative) $\eta$-directions a spiral $\vec{v}_{\pm}$ without self intersection, but with infinitely many points where $\vec{v}_+$ and $\vec{v_-}$ are equal (which all lie on the vertical axis\footnote{thanks to M. Hynek for this symmetry-argument, and sending me a picture of (part of) the $\gamma_0=-1$ curve},
\begin{figure}[h]
\centering
\includegraphics[width=8cm]{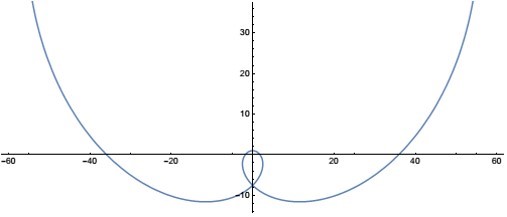}
\end{figure}\\
as $v_1(-\eta) = -v_1(\eta),\, v_2(-\eta) = v_2(\eta)$)
namely where $cos(\gamma_0 \eta)sinh(\eta) + \gamma_0 sin(\gamma_0 \eta)cosh(\eta) = 0$, i.e.
\begin{equation}\label{eq26}
\gamma_0 tan(\gamma_0 \eta) = -tanh(\eta).
\end{equation}
Alternatively:
$
R(\gamma_0 \eta_1)\big( \begin{smallmatrix}
s_1 \\ -\gamma_0 c_1
\end{smallmatrix}\big) =  R(\gamma_0 \eta_2)\big(\begin{smallmatrix}
s_2 \\ -\gamma_0 c_2
\end{smallmatrix} \big)
$
implies (considering the lengths of the 2 vectors) $\eta_2 = \pm \eta_1$, and $R(2\gamma_0 \eta)\big( \begin{smallmatrix} s \\ -\gamma_0 c \end{smallmatrix} \big) =
\big( \begin{smallmatrix} -s \\ -\gamma_0 c \end{smallmatrix} \big)$ then implies $\big( \begin{smallmatrix} \gamma_0 s_2 & c_2+1 \\ \gamma_0 (1-c_2) & s_2 \end{smallmatrix} \big)
\big(\begin{smallmatrix} c \\ s\end{smallmatrix} \big)
= \vec{0}
$, each of the 2 equations giving (\ref{eq26}).
Let us close with a few remarks: firstly, notice that (\ref{eq13}) implies (for $\gamma_0 \neq 0$) that the radius of curvature $\rho := \frac{1}{\kappa}$ and the (extended resp. signed) arclength parameter $s$ lie on a hyperbola:
\begin{equation}\label{eq27}
\kappa_0^2 \rho^2 - \mu^2 s^2 = \dfrac{\rho^2}{p^2} - \dfrac{s^2}{q^2} = \dfrac{\rho^2}{(\frac{1}{\gamma_0\omega\beta})^2} - \dfrac{s^2}{(\frac{1}{\omega\beta})^2} = 1,
\end{equation}
from which one can (see. e.g. \cite{6} immediately see that the curves (\ref{eq25}) are hypercycloids (studied as early as 1750, by Euler), arising by ``rolling a circle of complex radius $\zeta = \pm pq\frac{(p\pm iq)}{p^2+q^2}$ around one of purely imaginary radius $R = -i\frac{pq^2}{p^2+q^2}$'' \cite{6}.
Secondly, noting that (\ref{eq2}) and
\begin{equation}\label{eq28}
\begin{split}
\ddot{x} = \omega^2\begin{pmatrix}
-R\vec{u} \\ 0
\end{pmatrix},& \;x'' = \begin{pmatrix}
R\vec{u}'' \\ 0
\end{pmatrix},\\
\dot{x}' =  \omega^2\begin{pmatrix}
AR\vec{u}' \\ 0
\end{pmatrix}, & \;n \parallel \begin{pmatrix}
AR\vec{u}' \\ -\vec{u}\cdot \vec{u}' \omega
\end{pmatrix},
\end{split}
\end{equation}
i.e.
\begin{equation}\label{eq29}
h_{ab} \sim \begin{pmatrix}
\omega^2\vec{u} \times \vec{u}' & \omega \vec{u}'^2 \\
\cdot & \vec{u}'\times\vec{u}''
\end{pmatrix}
\end{equation}
imply that the vanishing of the mean curvature, $g^{ab}h_{ab} = 0$, reads
\begin{equation}\label{eq30}
\omega^2 \vec{u}'^2(\vec{u}\times \vec{u}') = (1+r^2\omega^2)\vec{u}'\times\vec{u}'',
\end{equation}
which is somewhat more complicated than the shape-equation, though of course equivalent: in arc-length  parametrization, $\vec{u}'^2 = 1$ one has $\vec{u}'\cdot \vec{u}'' = 0$, so
\begin{equation}\label{eq31}
\omega^2r sin(\phi) = (1+r^2\omega^2)|\vec{u}''| = (1+r^2\omega^2)\kappa,
\end{equation}
which, using (\ref{eq5}), (\ref{eq9}), (\ref{eq13}) and (\ref{eq16}) is easy to verify. To see the reverse is slightly more complicated (but of course true, as $H=0$ is equivalent to $\bigtriangleup x = 0$, hence in particular implies $\bigtriangleup t = 0$; note however the somewhat odd unbalance in simplicity concerning the 2 directions):
using $sin(\phi) = r\theta' = \sqrt{1-r'^2}$, and the lhs. of (\ref{eq13}) (which hold for \textit{any} arclength parametrized curve) one finds the second order ODE
\begin{equation}\label{eq32}
\dfrac{\omega^2 r^2}{1+\omega^2r^2}(1-r'^2) + (r''r + r'^2) = 1
\end{equation}
(with $\omega$ of course removable via $\omega r = \tilde{r}(\omega s)$); with u = $\frac{1}{2}\tilde{r}^2$ (\ref{eq32}) reads $(1+2u)u'' = (1+u'^2)$, resp. $(\ln(1+u'^2))' = 2\frac{u''u'}{1+u'^2} = \frac{2u'}{1+2u} = (\ln(1+2u))'$, implying
\begin{equation}\label{eq33}
1+\omega^2r^2r'^2 = (1+u'^2) = c^2(1+2u) = c^2(1+\omega^2r^2),
\end{equation}
which says that the first term on the lhs. of (\ref{eq32}) (hence also the second one) is separately constant (which is the shape-equation).
Thirdly, let us derive the Weierstrass-data by noting that (\ref{eq20}) can , with $c := a +ib$ and $w := u+iv$, be written as
\begin{equation}\label{eq34}
\begin{split}
x(u,v) & = \begin{pmatrix}
\frac{1}{2}c\, sinh(u+iv) + \frac{1}{2}\bar{c}\,sinh(u-iv) \\[0.15cm]
\frac{1}{2i}c \,cosh(u+iv) - \frac{1}{2i}\bar{c}\,cosh(u-iv)\\[0.15cm]
av+bu
\end{pmatrix}\\
 & =Re\int \vec{\varphi}(w)dw
\end{split}
\end{equation}
where
\begin{equation}\label{eq35}
\vec{\varphi} = c \begin{pmatrix}
cosh(w) \\ -isinh(w) \\ -i
\end{pmatrix}
\end{equation}
satisfies $\vec{\varphi}^2 = 0$ (guaranteeing resp. confirming harmonicity and isothemality). Fourthly, let us observe the following peculiarity: the planar curves $\vec{v}(\eta)$ (cp.(\ref{eq25})) not only satisfy the (reparametrization-invariant) shape equation,
\begin{equation}\label{eq36}
(\gamma_0^2 + 1)(\vec{v}\times \vec{v}')^2 = \vec{v}'^2\gamma_0^2(1+\vec{v}^2),
\end{equation}
but also
\begin{equation}\label{eq37}
\vec{v}'^2 = (\gamma_0^2 + 1)(\vec{v}^2+ 1), \; (\vec{v} \times \vec{v}') = \gamma_0(\vec{v}^2 + 1)
\end{equation}
(which reflect the conformal constancy of the metric); curiously, (\ref{eq37}) (though being parametrization-dependent) not only implies (\ref{eq36}) but also the \textit{linear} second order ODE
\begin{equation}\label{eq38}
\vec{v}'' - 2\gamma_0A\vec{v}' - (\gamma_0^2 + 1)\vec{v} = 0,
\end{equation}
which  is easily seen by differentiating (\ref{eq37}):
\begin{equation}\label{eq39}
\begin{split}
\vec{v}(2\gamma_0\vec{v}' + A\vec{v}'') = 0\; & \Rightarrow \;2\gamma_0\vec{v}' + A\vec{v}'' = \mu A \vec{v} \\
\vec{v}'(-\vec{v}'' + (\gamma_0^2+1)\vec{v}) = 0\; & \Rightarrow\; -\vec{v}''+ (\gamma_0^2+1)\vec{v} = \nu A \vec{v}'
\end{split}
\end{equation}
(from which $\mu = \gamma_0^2 +1,\, \nu = -2\gamma_0$, hence (\ref{eq38}), immediately follows). The general solution of (\ref{eq38}),
\begin{equation}\label{eq40}
\vec{v} = \left( \alpha_+ \begin{pmatrix}
c_0 \\s_0
\end{pmatrix} + \beta_+ \begin{pmatrix}
-s_0 \\c_0
\end{pmatrix} \right)e^{\eta} + \left( \alpha_- \begin{pmatrix}
c_0 \\s_0
\end{pmatrix} + \beta_- \begin{pmatrix}
-s_0 \\c_0
\end{pmatrix} \right)e^{-\eta}
\end{equation}
where $c_0 = cos(\gamma_0 \eta)$ and $s_0 = sin(\gamma_0 \eta)$, either inserted into (\ref{eq36}), or into (\ref{eq37}), yields (\ref{eq25}), resp.
\begin{equation}\label{eq41}
\alpha_+ = \dfrac{1}{2}, \;
\alpha_- = -\dfrac{1}{2}, \;
\beta_+ = -\dfrac{1}{2}\gamma_0, \;
\beta_- = -\dfrac{1}{2}\gamma_0.
\end{equation}
Fifthly, out of curiosity, consider the question whether in (\ref{eq20}) the first 2 components could possibly (without destroying isothemality) be more general linear combinations of the 4 harmonic functions $sinh(u)cos(v)$, $sinh(u)sin(v)$, $cosh(u)cos(v)$ and $cosh(u)sin(v)$, i.e.
\begin{equation}\label{eq42}
x_{\alpha=1,2} = (cos(v), sin(v))X_{\alpha}\begin{pmatrix}
cosh(u) \\ sinh(u)
\end{pmatrix},
\end{equation}
with constant $2\times 2$ matrices $X_1 = X$ and $X_2 = Y$. Orthogonality alone ($x_ux_v + y_uy_v = -ab$) gives the condition that $\big( P := \big(\begin{smallmatrix}
0&1 \\ 1&0
\end{smallmatrix}\big)\big)$
\begin{equation}\label{eq43}
\begin{split}
Z_{ij,kl} & := \big[\dfrac{1}{4}\sum_{\alpha}(X_{\alpha}P)_{ik}(AX_\alpha)_{jl} + (i\leftrightarrow j)\big] + \big[ k \leftrightarrow l\big]\\
& \stackrel {!}{=} ab\delta_{ij} (\sigma_3)_{kl}
\end{split}
\end{equation}
(which is invariant under the 3 obvious symmetries that result from $v \rightarrow v+v_0, u \rightarrow u_0$ or constant rotations in the $xy$-plane, i.e. in particular  $X_{\alpha} \rightarrow X_{\alpha}e^{u_0P}$ and/or $X_{\alpha} \rightarrow e^{v_0A}X_{\alpha}$); $\underset{j=k}{\sum}$ of (\ref{eq43}) e.g. gives the condition
\begin{equation}\label{eq44}
\sum_{\alpha}\big(
X_{\alpha}PAX_{\alpha} + Tr(X_{\alpha}P)AX_{\alpha} + X_{\alpha}PTr(AX_{\alpha}) + AX_{\alpha}^2P\big) = 4ab\sigma_3.
\end{equation}
While one can of course check that $X_1 = X = \big( \begin{smallmatrix} 0 & a \\ -b & 0 \end{smallmatrix} \big)$,
$X_2 = Y = \big( \begin{smallmatrix} b & 0 \\ 0 & a \end{smallmatrix} \big)$, satisfy the above equations (in (\ref{eq44}), e.g. $X$ giving $4ab\sigma_3 + (a^2-b^2)\mathbf{1}$, $Y(b^2-a^2)\mathbf{1}$), hence also all their symmetry-transforms, it is a priori nontrivial that the over determined nonlinear equations (\ref{eq43}) \textit{are} solvable, and have no other solutions than (\ref{eq25}), symmetry-transformed.\\\\
Finally, let us derive the shape-equation for hypersurfaces in \textit{arbitrary} dimensions, i.e. for
\begin{equation}\label{eq46}
x(t, \varphi^1 \ldots \varphi^M) = \begin{pmatrix}
e^{\Omega t}\vec{u}(\varphi^1 \ldots \varphi^M) \\ t
\end{pmatrix}
\end{equation}
with $\Omega$ a constant antisymmetric matrix. One finds
\begin{equation}\label{eq47}
\begin{split}
G_{AB} & = \begin{pmatrix}
\partial_a\vec{u}\partial_b\vec{u} & -\vec{u}\Omega\partial_a\vec{u} \\
\cdot & 1-\vec{u}\Omega^2\vec{u}
\end{pmatrix}  = \begin{pmatrix}
g_{ab} & u_a \\ u_b & \dot{x}^2
\end{pmatrix} \\
G & = g(1-\vec{u}\Omega^2\vec{u} - u_a g^{ab}u_b) =: gp^2 \\
G^{AB} & = \dfrac{1}{p^2}\begin{pmatrix}
p^2g^{ab} + u^au^b & -u^a \\ -u^b & 1
\end{pmatrix}\\
\sqrt{G}G^{AB} & = \dfrac{\sqrt{g}}{p}\begin{pmatrix}
p^2g^{ab} + u^au^b & -u^a \\ -u^b & 1
\end{pmatrix}
\end{split}
\end{equation}
where $u^a := g^{ac}u_c$. Letting $\bigtriangleup := \frac{1}{\sqrt{G}}\partial_A\sqrt{G}G^{AB}\partial_B$ act on $t$, the last component of $x$, one obtains the shape-equation
\begin{equation}\label{eq48}
\begin{split}
-\partial_a(\sqrt{g}\frac{g^{ab}}{p}u_b) & = 0 \\ \text{resp.} \qquad
+\vec{u}\Omega\partial_a(\sqrt{g}\frac{g^{ab}}{p}\partial_b\vec{u}) & = 0.
\end{split}
\end{equation}
Let us write $\vec{u} = r(\varphi^1 \ldots \varphi^M)\cdot \vec{e}(\varphi^1 \ldots \varphi^M)$ ($\vec{e}^2 = 1$), in analogy with $M=1$, where the solutions that we obtained could also be written as $\vec{u}(\theta) = r(\theta)\big( \begin{smallmatrix}
cos \theta \\ sin \theta \end{smallmatrix} \big)$
i.e. parametrized (instead of arclength) by the geometrical angle, for which (\ref{eq7}) gives, with $w = \omega^2r^2$ and for $\gamma_0(\theta-\theta_0)\geqslant 0$ (cp. II. 80 of \cite{3} for the relativistic analogue)
\begin{equation}\label{eq49}
\begin{split}
\theta-\theta_0 = \int d \theta & = \dfrac{1}{2} \gamma_0 \int \dfrac{dw}{w}\sqrt{\dfrac{w+1}{w-\gamma_0^2}} \\
 & = \gamma_0 (\gamma_0^2 + 1) \int \dfrac{d\psi}{\gamma_0^2 +(tanh \psi)^2} \\
 & = \gamma_0 \psi - arctan(\gamma_0 coth \psi)
 \end{split}
\end{equation}
upon $w = \gamma_0^2cosh^2 \psi + sinh^2 \psi$; so
\begin{equation}\label{eq50}
\begin{split}
tan(\theta-\theta_0) & = \dfrac{tan(\gamma_0\psi)-\gamma_0coth\psi}{1 + \gamma_0tan(\gamma_0\psi)coth\psi}\\
& = \dfrac{s_0sinh -\gamma_0c_0 cosh}{cos(\gamma_0 \psi)sinh \psi + \gamma_0 sin(\gamma_0 \psi)cosh \psi}
 \end{split}
\end{equation}
which (choosing $\theta_0 = \frac{\pi}{2}$, or $arccot$ instead of $-arctan$ in (\ref{eq49})) matches (cp. (\ref{eq21}))
\begin{equation}\label{eq51}
\dfrac{u_2}{u_1} = \dfrac{cos(\gamma_0 v)sinh (v) + \gamma_0sin(\gamma_0v)cosh(v)}{sin(\gamma_0 v)sinh (v) - \gamma_0cos(\gamma_0v)cosh(v)},
\end{equation}
$v = arcsinh(\mu s) = \psi$, i.e. $sinh \psi = \mu s$, confirming (cp.(\ref{eq9})) $\omega^2 r^2 = w = \gamma_0^2(1+\mu^2s^2) + \mu^2s^2 = \gamma_0^2 + \frac{\omega^2 s^2}{1+\gamma_0^2}$.\\\\
Consider now arbitrary $M \geqslant 1$: with
\begin{equation}\label{eq52}
\begin{split}
\vec{u} & = r(\varphi^1 \ldots \varphi^M)\vec{e}(\varphi^1 \ldots \varphi^M),\; \vec{e}^2 = 1 \\
g_{ab} & = \partial_a\vec{u}\partial_b\vec{u} = r^2(\partial_a\vec{e}\partial_b \vec{e} + \dfrac{r_a r_b }{r^2}) =: r^2(\bar{g}_{ab} + v_a v_b)\\
g^{ab} & = \dfrac{1}{r^2}(\bar{g}^{ab} - \dfrac{v^av^b}{q^2}),\; q := \sqrt{1 + v_av_b\bar{g}^{ab}} \\
g & = r^{2M}\bar{g}q^2
\end{split}
\end{equation}
the shape-equation, (\ref{eq48}), becomes
\begin{equation}\label{eq53}
\vec{e}\Omega \partial_a(\bar{f}\sqrt{\bar{g}}\underbrace{(q^2\bar{g}^{ab} - v^a v^b)}_{=: \rho^{ab} = q^2 f^{ab}}\partial_b\vec{e}) = 0
\end{equation}
with $\bar{f} = \frac{r^M}{qp} = \frac{r^{M+1}}{\sqrt{(r^2 + r^a r_a)[1+ (\Omega\vec{e})^2-e^2] + r^4(e \cdot v)^2}}$,
$(\Omega\vec{e})^2-e^2  =0 \,\forall \Omega = -\Omega^T$; indices are raised and lowered with $\bar{g}^{ab}$, resp. $\bar{g}_{ab}$ and $e^2 = e_a\bar{g}^{ab}e_b$, $e\cdot v = e_a\bar{g}^{ab}v_b$, $e_a := \vec{e}\Omega \partial_a \vec{e}$.
For $M = 1$, $\vec{e} = \big( \begin{smallmatrix} cos \varphi \\ sin \varphi \end{smallmatrix} \big)$, $\rho^{ab} = 1 = \sqrt{\bar{g}} = \bar{g}_{ab}$, $e_a = e_1 =\omega$, $\Omega = \omega \big( \begin{smallmatrix} 0 & -1 \\ 1 & 0 \end{smallmatrix} \big)$ and $(\bar{f})^2 =$ const. $= c^2$ gives $r^4 = c^2[(r^2 + r'^2) + \omega^2r'^2r^2]$ which upon $c^2 = \frac{\gamma_0^2}{\omega^2}$ gives
\begin{equation}\label{eq54}
r'^2 = \dfrac{r^2(\omega^2 r^2 - \gamma_0^2)}{\gamma_0^2(1+\omega^2r^2)},
\end{equation}
in agreement with (\ref{eq49}).\\\\
Before rewriting (\ref{eq53}), $\vec{e}$ parametrizing a unit $M$-dimensional sphere, let us note that the solutions of theorem 1 in \cite{7} are of the form (\ref{eq46}), hence should satisfy the $M = 2k+1$ dimensional shape-equation (\ref{eq48}) [resp. (\ref{eq53})] for the particular case of
\begin{equation}\label{eq55}
\begin{split}
r(\varphi^1 \ldots \varphi^{2k+1 = M}) & = \varphi^M, \\
\vec{e}  = \vec{m}(\varphi^1 \ldots \varphi^{2k})  & = \dfrac{1}{\sqrt{2}} \begin{pmatrix}
\vec{m}_1(\varphi^1 \ldots \varphi^k) \\ \vec{m}_2(\varphi^{k+1} \ldots \varphi^{2k}
\end{pmatrix}.
\end{split}
\end{equation}
As $\vec{e}$ in that case does not depend on one of the parameters, $(\varphi^M)$, resp.
\begin{equation}\label{eq56}
\begin{split}
g_{ab} & = \partial_a(r\vec{m})\partial_b(r\vec{m}) =  \begin{pmatrix}
\bar{g}_{\alpha \beta} = \partial_{\alpha}\vec{m}\partial_{\beta}\vec{m} & 0 \\ 0 & 1
\end{pmatrix}, \\
g^{ab} & = \begin{pmatrix}
\bar{g}^{\alpha \beta} & 0 \\ 0 & 1
\end{pmatrix},\; g = \bar{g} \\
u_M & = -r\vec{m} \Omega \partial_M(r\vec{m}) = 0
\end{split}
\end{equation}
(and $-r^2\vec{m} \Omega \partial_{\alpha}\vec{m} = -r^2 e_{\alpha}, \, \alpha = 1, \ldots. 2k$) (\ref{eq48}) reduces (in the case that $\vec{m}$ describes a minimal surface in $S^{2k}$, i.e. satisfies $\frac{1}{\sqrt{\bar{g}}}\partial_{\alpha}\sqrt{\bar{g}}\bar{g}^{\alpha \beta}\partial_{\beta}\vec{m} \equiv 2k\vec{m}$) to the interesting condition\footnote{Thanks to J.Arnlind for a discussion on the geometric meaning,
as well as to M.Bordemann, F.Finster, and G.Linardopoulos for more general discussions}
\begin{equation}\label{eq57}
e^{\alpha}\partial_{\alpha}e^2 = 0
\end{equation}
where $e^{\alpha} = \bar{g}^{\alpha \beta}e_{\beta}$ and $e^2 = e_{\alpha}\bar{g}^{\alpha \beta}e_{\beta}$; the proof of theorem 1 of \cite{7} can thus be much simplified as in order to verify (\ref{eq57}) one does not need to calculate the second fundamental form at all, while only needing the technical Lemma \cite{7}
\begin{equation}\label{eq58}
\begin{split}
\dfrac{1}{2}(\vec{m}_1\partial_{\alpha}\vec{m}_2)\bar{g}^{\alpha \beta}(m_1 \partial_{\alpha}\vec{m}_2)
 & = 1-(\vec{m}_1 \cdot \vec{m}_2)^2 \\
 & = \dfrac{1}{2}(\vec{m}_2\partial_{\alpha}\vec{m}_1)\bar{g}^{\alpha \beta}(m_2 \partial_{\beta}\vec{m}_1)
\end{split}
\end{equation}
to see that (using also the block diagonality of $\bar{g}^{\alpha \beta}$)
\begin{equation}\label{eq59}
e^2  = e_{\alpha}\bar{g}^{\alpha \beta} e_{\beta} = 1 - (\vec{m}_1 \cdot \vec{m}_2)^2,
\end{equation}
and (using again (\ref{eq58}))
\begin{equation}
\begin{split}
e^{\alpha}\partial_{\alpha}e^2 & = -2(\vec{m}_1 \cdot \vec{m}_2)\cdot
(\vec{m}_2 \partial_{\beta} \vec{m}_1 - \vec{m}_1 \partial_{\beta} \vec{m}_2)\bar{g}^{\alpha \beta} \cdot
(\vec{m}_2 \partial_{\alpha} \vec{m}_1 + \vec{m}_1 \partial_{\alpha} \vec{m}_2) \\
 &
= -4(\vec{m}_1 \cdot \vec{m}_2)\{ (1- (\vec{m}_1 \cdot \vec{m}_2 )^2) - (1- (\vec{m}_1 \cdot \vec{m}_2 )^2) \} = 0.
\end{split}
\end{equation}
Note that (\ref{eq57}) can be simplified as
$$
\partial_{\gamma}\vec{m}\nabla_{\alpha \beta}^2 \vec{m} = \partial{\gamma}\vec{m}(\partial_{\alpha \beta}^2\vec{m} - \gamma_{\alpha \beta}\partial_{\varepsilon}\vec{m}) = 0
$$
implies
\begin{equation}\label{eq60}
\begin{split}
\dfrac{1}{2} e ^{\alpha}\nabla_{\alpha}(e^{\beta}e_{\beta}) & = e^{\alpha}e^{\beta}\nabla_{\alpha}e_{\beta} \\
 & = e^{\alpha}e^{\beta} [\partial_{\alpha}\vec{m}\Omega\partial_{\beta}\vec{m} + \vec{m}\Omega \nabla_{\alpha}\nabla_{\beta}m]\\
 & = e^{\alpha}e^{\beta} [0 - e^{\gamma}\partial_{\gamma}\vec{m} \nabla_{\alpha}\nabla_{\beta}\vec{m} + \rho\vec{n}\nabla_{\alpha}\nabla_{\beta}\vec{m}] \\
 & = \rho e^{\alpha}e^{\beta} \tilde{h}_{\alpha \beta} \stackrel{!}{=} 0,
\end{split}
\end{equation}
meaning that $\Omega$ and the minimal co-dimension 1 hypersurface $\vec{m} \in S^{2k+1}$ have to be such that $\Omega\vec{m} = -e^{\gamma}\partial_{\gamma}\vec{m} + \rho \vec{n}$ is purely tangential (i.e. $\rho = 0$), \textit{or} $e^{\alpha}$ is a null direction for the second fundamental form $\tilde{h}_{\alpha \beta} := \vec{n}\gamma_{\alpha \beta}^2 \vec{e}$ $(\vec{n}\cdot \vec{m} = 0 = \vec{n} \partial_{\alpha}\vec{m})$.\\
Naturally one would like to know for which $\vec{m}$ and $\Omega$ the condition (\ref{eq60}) holds; let us consider
\begin{equation}\label{eq61}
\vec{m} = \dfrac{1}{\sqrt{p+q}}\begin{pmatrix}
\sqrt{q}\vec{m}_1(\theta_1 \ldots \theta_p) \\
\sqrt{p}\vec{m}_2(\varphi_1 \ldots \varphi_q)
\end{pmatrix},
\end{equation}
with $\vec{m}_1$ and $\vec{m}_2$ parametrizations of the unit spheres $S^p$ resp. $S^q$ (making $\vec{m}$ a minimal hypersurface in $S^{p+q=1}$), $
\Omega = \big( \begin{smallmatrix} \Omega' & B \\ -B^T & \Omega''
\end{smallmatrix} \big)
= -\Omega^T
$, with $B$ an arbitrary constant $(p+1)\times (q+1)$ matrix. $\rho$ (cp. (\ref{eq60})) can be zero only if $B=0$, as $\vec{n} =  \frac{1}{\sqrt{p+q}}
\big( \begin{smallmatrix}\sqrt{p}\vec{m}_1 \\ -\sqrt{q}\vec{m}_2
\end{smallmatrix} \big)
$ and
\begin{equation}\label{eq62}
\Omega \vec{m} = \dfrac{1}{\sqrt{p+q}}\begin{pmatrix}
\sqrt{q}\Omega'\vec{m}_1 + \sqrt{p}B\vec{m}_2 \\
-\sqrt{q}B^T\vec{m}_1 + \sqrt{p}\Omega''\vec{m}_2
\end{pmatrix}
\end{equation}
gives
\begin{equation}\label{eq63}
(p+q)\vec{n}\Omega \vec{m} = p\vec{m}_1 B \vec{m}_2 + q\vec{m}_2 B^T \vec{m}_1
= (p+q)(\vec{m}_1 B \vec{m}_2).
\end{equation}
On the other hand,
\begin{equation}\label{eq64}
\begin{split}
(p+q)e_{\alpha} & = q \vec{m}_1 \Omega' \partial_{\alpha} \vec{m}_1 + \sqrt{pq}(\vec{m}_1 B \partial_{\alpha} \vec{m}_2 - \vec{m}_2 B^T \partial_{\alpha} \vec{m}_1) \\
& \quad + p \vec{m}_2 \Omega'' \partial_{\alpha} \vec{m}_2, \\
g_{\alpha \beta} & = \dfrac{1}{p+q} \begin{pmatrix}
q\bar{g}_{\alpha' \beta'} & 0 \\ 0 & p\bar{g}_{\alpha'' \beta''}
\end{pmatrix},\\
\tilde{h}_{\alpha \beta} & = \vec{n}\partial_{\alpha \beta}^2 \vec{m} = \dfrac{1}{p+q}\{ \sqrt{pq}\vec{m}_1 \partial_{\alpha \beta}^2 \vec{m}_1 - \sqrt{pq}\vec{m}_2 \partial_{\alpha \beta}^2 \vec{m}_2\}\\
 & =  \dfrac{\sqrt{pq}}{p+q} \begin{pmatrix}
-\bar{g}_{\alpha' \beta'} & 0 \\ 0 & \bar{g}_{\alpha'' \beta''}
\end{pmatrix}
\end{split}
\end{equation}
imply
\begin{equation}\label{eq65}
\begin{split}
e^{\alpha'} & = g^{\alpha' \beta} e_{\beta} = \bar{g}^{\alpha' \beta'}(\vec{m}_1 \Omega' \partial_{\beta'} \vec{m}_1 -
\sqrt{\frac{p}{q}}\vec{m}_2 B^T \partial_{\beta'} \vec{m}_1)\\
e^{\alpha''} & = g^{\alpha'' \beta} e_{\beta} = \bar{g}^{\alpha'' \beta''}(\vec{m}_2 \Omega'' \partial_{\beta''} \vec{m}_2 +
\sqrt{\frac{q}{p}}\vec{m}_1 B \partial_{\beta''} \vec{m}_2),
\end{split}
\end{equation}
hence $e^{\alpha}e^{\beta}\tilde{h}_{\alpha \beta} = 0$ equivalent to $e^{\alpha'}e^{\beta'}\bar{g}_{\alpha' \beta'} =: e'^2
\stackrel{!}{=} e''^2 := e^{\alpha''}e^{\beta''}\bar{g}_{\alpha'' \beta''}
$, resp.
\begin{equation}\label{eq66}
\begin{split}
\bar{g}^{\alpha' \beta'} (\vec{m}_1 \Omega' \partial_{\alpha'} \vec{m}_1 -
\sqrt{\frac{p}{q}}\vec{m}_2 B^T \partial_{\alpha'} \vec{m}_1)
(\vec{m}_1 \Omega' \partial_{\beta'} \vec{m}_1 -
\sqrt{\frac{p}{q}}\vec{m}_2 B^T \partial_{\beta} \vec{m}_1)\\
\stackrel{!}{=}
\bar{g}^{\alpha'' \beta''} (\vec{m}_2 \Omega'' \partial_{\alpha''} \vec{m}_2 +
\sqrt{\frac{q}{p}}\vec{m}_1 B \partial_{\alpha''} \vec{m}_2)
(\vec{m}_2 \Omega'' \partial_{\beta''} \vec{m}_2 +
\sqrt{\frac{q}{p}}\vec{m}_1 B \partial_{\beta''} \vec{m}_2).
\end{split}
\end{equation}
Putting $\Omega' = 0 = \Omega''$ (they only generate geometrically insignificant rotations of the 2 individual spheres) one gets the condition
\begin{equation}\label{eq67}
q\bar{g}^{\alpha'' \beta''} (\vec{m}_1 B \partial_{\alpha''} \vec{m}_2)
(\vec{m}_1 B \partial_{\beta''} \vec{m}_2)
\stackrel{!}{=}
p\bar{g}^{\alpha' \beta'} (\vec{m}_2 B^T \partial_{\alpha'} \vec{m}_1)
(\vec{m}_2 B^T \partial_{\beta'} \vec{m}_1)
\end{equation}
which can be reformulated as : The projection of $\sqrt{p}B\vec{m}_2 \in \mathbb{R}^{p+1}$ onto $T_{\theta}S^p$ (the tangent space of $S^p$ at $\theta$) has to have the same length then $\sqrt{q}B^T\vec{m}_1 \in \mathbb{R}^{q+1}$ projected onto $T_{\varphi}S^q$ (with $Q := \bar{g}^{\alpha''\beta''}\partial_{\alpha''}\vec{m}_2\partial_{\beta''}\vec{m}_2^T$, $Q^T = Q$, $Q^2 = Q$ just as $P := \bar{g}^{\alpha'\beta'}\partial_{\alpha'}\vec{m}_1\partial_{\beta'}\vec{m}_1^T$, $P^T = P$, $P^2 = P$ projects onto $TS^p$); equivalently, with $\mathbb{Q}_2 = qBQB^T = q(BQ)(BQ)^T$, $\mathbb{Q}_1 = pB^TPB = p(PB)^T(PB)$, $\vec{m}_1^T\mathbb{Q}_2\vec{m}_1 \stackrel{!}{=} \vec{m}_2 \mathbb{Q}_1 \vec{m}_2$.\\
Expanding $B\vec{m}_2$ and $B^T\vec{m}_1$ into their tangent and normal parts,
$$
B\vec{m}_2 = (\vec{m}_1 B \vec{m}_2)\vec{m}_1 + \bar{g}^{\alpha' \beta'}(\partial_{\alpha'}\vec{m}_1 B \vec{m}_2)(\partial_{\beta'}\vec{m}_1 B \vec{m}_2)
$$
$$
B^T\vec{m}_1  = (\vec{m}_2 B^T \vec{m}_1)\vec{m}_2 + \bar{g}^{\alpha'' \beta''}(\partial_{\alpha''}\vec{m}_2 B^T \vec{m}_1)(\partial_{\beta''}\vec{m}_2 B^T \vec{m}_1),
$$
and squaring, one obtains the condition
\begin{equation}\label{eq68}
\begin{split}
q(\nu - (\vec{m}_2 B^T \vec{m}_1)^2) \stackrel{!}{=} p(\mu - (\vec{m}_1 B \vec{m}_2)^2)\\
\nu := \vec{m}_1^T B B^T \vec{m}_1 = \nu(\theta_1 \ldots \theta_p), \quad 
\mu := \vec{m}_2 B^T B \vec{m}_2 = \mu(\varphi_1 \ldots \varphi_q)
\end{split}
\end{equation}
which for $p = q$ is indeed satisfied if $B$ is a multiple of an orthogonal matrix (like e.g. $B = -\mathbf{1}$, as in \cite{7}), but\footnote{responding to a question J. Choe on several occasions asked} seems impossible to satisfy for $p\neq q$.\\
Can one find explicit non-trivial solutions of (\ref{eq53}) for $M > 1$ ? As $\vec{e}$, as a parametrization of $S^M$, satisfies
\begin{equation}\label{eq69}
\underline{\vartriangle} \vec{e} = \dfrac{1}{\sqrt{\bar{g}}}\partial_a(\sqrt{\bar{g}}\bar{g}^{ab}\partial_b \vec{e})
= -M\vec{e}\, \bot \,\Omega \vec{e}
\end{equation}
(\ref{eq53}) can be written as
\begin{equation}\label{eq70}
\begin{split}
0 & = f^{ab} e_b \{\partial_a \ln \bar{f} - [(\bar{g}_{ac} + v_a v_c) \bigtriangledown_e (\dfrac{v^e v^c}{q^2}) - \partial_a \ln q^2] \} - \dfrac{v^a v^c}{q^2}\bigtriangledown_a e_c\\
 & = f^{ab}e_b\{\partial_a \ln \bar{f} + \dfrac{1}{2}\partial_a \ln q^2 -v_a [\underbrace{\bigtriangledown_e v^e - \dfrac{v^e v^c}{q^2}\bigtriangledown_e v_c}_{f^{ec}\bigtriangledown_e \bigtriangledown_c \ln r =: \hat{\Delta} \ln r} ] \}- \dfrac{v^a v^c}{q^2}\bigtriangledown_a e_c
\end{split}
\end{equation}
with $\bar{f}^2 = \frac{r^{2M}}{q^2 + r^2(e\cdot v)^2} =: \frac{1}{Y^2}$, i.e
\begin{equation}\label{eq71}
q^2 + r^2(e \cdot v)^2 = r^{2M}Y^2.
\end{equation}
(\ref{eq70}), resp.
\begin{equation}\label{eq72}
\dfrac{e \cdot v}{q^2} \hat{\Delta} \ln r = \dfrac{1}{2}f^{ab}e_b \partial_a \ln(\bar{f}^2 q^2)- \dfrac{v^a v^c}{q^2}\bigtriangledown_a e_c
\end{equation}
becomes ( if $e \cdot v \neq 0$)
\begin{equation}\label{eq73}
\Delta \ln r := (\bigtriangledown^a\bigtriangledown_a \ln r) = \dfrac{1}{2}(\dfrac{e^a \bigtriangledown_a}{e \cdot v}q^2) -
\dfrac{1}{2}\dfrac{q^2 f^{ab}}{e \cdot v} e_b \partial_a \ln Y^2 - \dfrac{v^a v^c}{e \cdot v}\bigtriangledown_a e_c,
\end{equation}
the last but one term being $-\frac{rY'(r)}{Y}$ if $Y = Y(r)$ (and $e\cdot v \neq 0$); and
\begin{equation}\label{eq74}
\begin{split}
\Delta \ln r & = \dfrac{e^a \bigtriangledown_a}{e \cdot v}q - \dfrac{v^a v^c}{e \cdot v}\bigtriangledown_a e_c = \dfrac{v^a}{e \cdot v}\bigtriangledown_a (e \cdot v) \\
\text{i.e.} \quad \bigtriangledown_a \big( \dfrac{v^a}{e \cdot v}\big) & = 0
\end{split}
\end{equation}
if $Y = const$ (and $e\cdot v \neq 0$). For $M=1$, (\ref{eq74}) is trivially satisfied (and in fact \textit{all} solutions of the shape-equation are given by (\ref{eq71}), with $Y = const$).
The case $e \cdot v = 0$, which we excluded in the last 2 equations
is the one special case when the shape equation does \textit{not} imply that all the other components of $\Delta x = 0$ are satisfied, for the following reason: For \textit{any} $x(t, \varphi^1 \ldots \varphi^M)$ (not necessarily minimal) one has, for $x$ of the form (\ref{eq46}),$\dot{x}\cdot \Delta x = 0 = \partial_a x \cdot \Delta x$;
\begin{equation}\label{eq75}
\begin{split}
\text{the first equation gives } \qquad & \Delta t + (\Omega R \vec{u})^T \Delta R\vec{u} = 0 \\
\text{and the second} \qquad \qquad \qquad & (R\partial_a \vec{u})^T \Delta R \vec{u}  = 0.
\end{split}
\end{equation}
As long as the $N = M + 1$ vectors $\Omega \vec{u}$ and $\partial_a \vec{u} = \vec{e}\partial_a r + r \partial_a \vec{e}$ are linearly independent, $\Delta t $ will therefore imply $\Delta R \vec{u} = \vec{0}$. Due to $\Omega \vec e = - e^a\partial_a \vec{e}$
the equation $\vec{0} = \mu^a \partial_a \vec{u} + \mu r \Omega \vec{e} = r(\mu_a - e^a \mu)\partial_a \vec{e} + (\mu^a \partial_a r)\vec{e}$ can have a non-trivial solution, namely $\mu^a = \mu e^a$, if $e \cdot v = 0$. If on the other hand $e \cdot v  \neq 0$, $\Delta t = 0$ (the shape-equation) clearly implies $\Delta x = 0$. Let us nevertheless consider (explicitly) an example with $e \cdot v = 0$; $N = 4$, $\Omega = \big( \begin{smallmatrix} A & 0 \\ 0 & A \end{smallmatrix} \big)$, now
$A  = \big( \begin{smallmatrix} 0 & -1 \\ 1 & 0 \end{smallmatrix} \big)$
\begin{equation}\label{eq76}
\begin{split}
\vec{e} & = \begin{pmatrix}
sin \theta_1\, sin \theta_2 \, cos \theta_3 \\
sin \theta_1\, sin \theta_2 \, sin \theta_3 \\
sin \theta_1\, cos \theta_2 \\
cos \theta_1
\end{pmatrix}, \quad
\bar{g}_{ab} = \begin{pmatrix}
1 & 0 & 0 \\
0 & sin^2 \theta_1 & 0 \\
0 & 0 & sin^2\theta_1\, sin^2\theta_2
\end{pmatrix} \\[0.10 cm]
e_1 & = -cos\theta_2, \; e_2 = sin \theta_1\, cos \theta_1 \, sin \theta_2, \; e_3 = sin^2 \theta_1\, sin^2 \theta_2 \\
e^1 & = -cos\theta_2, \; e^2 = cot \theta_1\, sin \theta_2, \; e^3 = 1.
\end{split}
\end{equation}
A $\theta_3$ independent solution of $e^a \partial_a r = 0$ is
\begin{equation}\label{eq77}
r(\theta_1\theta_2\theta_3) = d^2(sin\theta_1 sin\theta_2)^k
\end{equation}
($k$ being the separation constant resulting from $(-cos\theta_2 \partial_1 + cot\theta_1sin\theta_2\partial_2)f_1(\theta_1)f_2(\theta_2) = 0$ resp. $tan \theta_i \frac{f_i'}{f_i} = \lambda$). What about the other special case, $(v^a) \parallel (e^a)$? \\
To include/cover the case $e \cdot v = 0$, let us derive the mean curvature  $H$ of (\ref{eq46}) (as interested in $H = 0$, not caring about overall constants, related e.g. to the normalization of the hypersurface normal, $n$)
\begin{equation}\label{eq78}
\begin{split}
\dot{x} & = \begin{pmatrix}
\Omega R \vec{u} \\ 1
\end{pmatrix}, \qquad \partial_a x = \begin{pmatrix}
R \partial_a \vec{u} \\ 0
\end{pmatrix}, \\
n & \parallel \begin{pmatrix}
R \vec{n} := \frac{R(\vec{e} - v^a \partial_a \vec{e})}{\sqrt{1+v^2}} \\[0.15 cm]
z = -\frac{r}{q}(e \cdot v) = -\nu
\end{pmatrix},
\end{split}
\end{equation}
where $z$ is not needed (as $\partial_{AB}^2 t = 0$), and
\begin{equation}\label{eq79}
\vec{n} \partial_a \vec{u} = \dfrac{(\vec{e} - v^c \partial_c \vec{e})}{\sqrt{1+v^2}}(r\partial_a \vec{e} + \vec{e}\partial_a r) = 0,
\end{equation}
i.e. $\vec{n}$ normal to the tangent plane of the hypersurface in $\mathbb{R}^{N = M+1}$ described by $\vec{u}(\varphi^1 \ldots \varphi^M)$. Using (\ref{eq47}) the minimality of (\ref{eq46}) reads
\begin{equation}\label{eq80}
\vec{n} \Omega^2 \vec{u} - 2 u^a \vec{n} \Omega \partial_a \vec{u} + (p^2g^{ab} + u^a u^b) \vec{n} \partial_{ab}^2 \vec{u} = 0
\end{equation}
where
$
u^a := - g^{ab}\vec{u}\Omega \partial_a \vec{u} \underset{\text{if}(e \cdot v = 0)}{=} - e^a := -\bar{g}^{ab}(\vec{e}\Omega \partial_b \vec{e})
$. \\
From now on always assuming $e \cdot v = 0$, one finds
\begin{equation}\label{eq81}
q \cdot \vec{n} \partial_{ab}^2 \vec{u} = \bigtriangledown_a \bigtriangledown_b r - 2 v_a v_b r - r \bar{g}_{ab}, \quad p^2 =1,
\end{equation}
so that the third term in (\ref{eq80}) equals
\begin{equation}\label{eq82}
\dfrac{1}{r^2}\left(\bigtriangleup r - M r - \dfrac{v^2 r}{q^2} - \dfrac{v^av^b}{q^2}\bigtriangledown_a \bigtriangledown_b r\right)
+ e^a e^b \bigtriangledown_a \bigtriangledown_b r - r e^2,
\end{equation}
while
\begin{equation}\label{eq83}
\begin{split}
q\vec{n} \Omega^2 r \vec{e} & = -e(\Omega r \vec{e})^2 r - r v^a \vec{e} \Omega^2 \partial_a \vec{e} \\
2e^a(\vec{e} - v^c \partial_c \vec{e})\Omega\partial_a \vec{u} & = 2r \{ e^2 - e^av^c \partial_c \vec{e}\Omega \partial_a \vec{e} \}.
\end{split}
\end{equation}
Using $\Omega \vec{e} = -e^c \partial_c \vec{e}$, $e_b = \vec{e} \Omega \partial_b \vec{e}$
\begin{equation}\label{eq84}
\partial_a \vec{e} \Omega^2 \vec{e} = -e^c \partial_a \vec{e} \Omega \partial_c \vec{e},
\end{equation}
hence altogether 
\begin{equation}\label{eq85}
\bigtriangleup r -M r - \dfrac{v^2 r}{q^2} - \dfrac{v^a v^b}{q^2}\bigtriangledown_a \bigtriangledown_b r + r^2 e^a e^b \bigtriangledown_a \bigtriangledown_b r + r^3 v^a \partial_a \vec{e}\Omega^2 \vec{e} = 0
\end{equation}
the last term being zero if $\Omega^2 = -\omega^2 \mathbf{1}$, and (always) the last 2 terms combining to give
$-r^3 e^a v^b (\vec{e} \Omega \bigtriangledown_a \bigtriangledown_b \vec{e}) \stackrel{?}{=} 0$. Multiplying by $r$ and integrating, the first 3 terms are manifestly negative, and $r^3 e^a e^b \bigtriangledown_a (r v_b) = -r^4 e^a v^b \bigtriangledown_a e_b$. Note also that $\bigtriangleup \vec{e} := \frac{1}{\sqrt{\bar{g}}}\partial_a\sqrt{\bar{g}}\bar{g}^{ab}\partial_b \vec{e} (= \bigtriangledown^a \bigtriangledown_a \vec{e}) = -M\vec{e}$, projected onto $\Omega \vec{e}$, implies
\begin{equation}\label{eq86}
\bigtriangledown^a e_a = \bar{g}^{ab} \partial_a \vec{e} \Omega \partial_b \vec{e} = 0.
\end{equation}
How about $M = 2$ ? (in (\ref{eq76})/(\ref{eq77})) we considered $M = 3$, as for even $M$ $\Omega$ will always have a zero eigenvalue); taking $\Omega = \omega \left( \begin{smallmatrix}
0 & -1 & 0 \\
1 & 0 & 0 \\
0 & 0& 0
\end{smallmatrix} \right)
$, $\vec{e} = \left( \begin{smallmatrix}
sin \theta\,cos\varphi \\
sin \theta\,sin\varphi \\
cos \theta
\end{smallmatrix} \right) $
one gets $(\bar{g}_{ab}) = \left( \begin{smallmatrix}
1 \\ sin^2 \theta
\end{smallmatrix} \right)$,
$\Omega \vec{e} = \partial_{\varphi} \vec{e}$; $e_1 = 0$, $e_2 = -sin^2 \theta$, $e^1 = 0$, $e^2 = -1$; so any $\varphi$ independent $r$ will have $v \cdot e = 0$. The effect of $e^{\Omega t}$ corresponds to $\varphi \rightarrow \varphi + t\omega$, i.e. can be absorbed by a reparametrization; hence $x(\tilde{\varphi}, \theta, t) =  \left( \begin{smallmatrix}
\underset{^{\sim}}{\vec{u}} (\theta, \tilde{\varphi}) \\t
\end{smallmatrix} \right)$, i.e. giving only a trivial direct product hypersurface.
Finally, note that the above also applies to timelike minimal hypersurfaces in Minkowskispace, if the following changes are made: replace $p$ by $\tilde{p}:= \sqrt{1 - r^2\frac{(e \cdot v)^2}{q^2}}$, (\ref{eq47}) by
$(\sqrt{G}G^{\alpha \beta})_{\alpha\beta = 0\ldots M} =
\frac{\sqrt{g}}{\tilde{p}}
\left( \begin{smallmatrix}
1 & -u^b \\
-u^a & u^a u^b - \tilde{p}^2 g^{ab}
\end{smallmatrix} \right);
$ correspondingly: in (\ref{eq53}), and (\ref{eq70})-(\ref{eq73}):
$
\bar{f} \rightarrow \bar{f}_{rel} = \frac{r^M}{\tilde{p}q} = \frac{1}{\tilde{Y}}
$,
and in (\ref{eq80}): replace $p^2$ by $-\tilde{p}^2$, hence in (\ref{eq82}) and (\ref{eq85}) all terms that do not contain $\Omega$ (resp.$e's$) changing sign.

\end{document}